\documentclass[twoside]{article}
\include{graphicx}
\title{The asymptotics of the Touchard polynomials}
\author{\sc R. B.\ Paris \\
{\em Division of Computing and Mathematics}, \\
{\em University of Abertay Dundee, Dundee DD1 1HG, UK}
}
\begin{document}
\def\f#1#2{\mbox{${\textstyle \frac{#1}{#2}}$}}
\def\dfrac#1#2{\displaystyle{\frac{#1}{#2}}}
\def\boldal{\mbox{\boldmath $\alpha$}}
{\newcommand{\Sgoth}{S\;\!\!\!\!\!/}
\newcommand{\bee}{\begin{equation}}
\newcommand{\ee}{\end{equation}}
\newcommand{\lam}{\lambda}
\newcommand{\ka}{\kappa}
\newcommand{\al}{\alpha}
\newcommand{\fr}{\frac{1}{2}}
\newcommand{\fs}{\f{1}{2}}
\newcommand{\g}{\Gamma}
\newcommand{\br}{\biggr}
\newcommand{\bl}{\biggl}
\newcommand{\ra}{\rightarrow}
\newcommand{\mbint}{\frac{1}{2\pi i}\int_{c-\infty i}^{c+\infty i}}
\newcommand{\mbcint}{\frac{1}{2\pi i}\int_C}
\newcommand{\mboint}{\frac{1}{2\pi i}\int_{-\infty i}^{\infty i}}
\newcommand{\gtwid}{\raisebox{-.8ex}{\mbox{$\stackrel{\textstyle >}{\sim}$}}}
\newcommand{\ltwid}{\raisebox{-.8ex}{\mbox{$\stackrel{\textstyle <}{\sim}$}}}
\renewcommand{\topfraction}{0.9}
\renewcommand{\bottomfraction}{0.9}
\renewcommand{\textfraction}{0.05}
\newcommand{\mcol}{\multicolumn}
\date{}
\maketitle
\pagestyle{myheadings}
\markboth{\hfill \sc R. B.\ Paris  \hfill}
{\hfill \sc  Touchard polynomials\hfill}
\begin{abstract}
We examine the asymptotic expansion of the Touchard polynomials $T_n(z)$ (also known as the exponential polynomials)
for large $n$ and complex values of the variable $z$. In our treatment $|z|$ may be finite or allowed to be large like $O(n)$. We employ the method of steepest descents to a suitable integral representation of $T_n(z)$ and find that the number of saddle points that contribute to the expansion depends on the values of $n$ and $z$.
Numerical results are given to illustrate the accuracy of the various expansions.

\vspace{0.4cm}

\vspace{0.3cm}

\end{abstract}

\vspace{0.3cm}

\noindent $\,$\hrulefill $\,$

\vspace{0.2cm}

\begin{center}
{\bf 1. \  Introduction}
\end{center}
\setcounter{section}{1}
\setcounter{equation}{0}
\renewcommand{\theequation}{\arabic{section}.\arabic{equation}}
The Touchard polynomials $T_n(z)$, also known as exponential polynomials, are defined by
\bee\label{e11}
T_n(z)=e^{-z}\sum_{k=0}^\infty \frac{k^nz^k}{k!}=e^{-z} \bl(z \frac{d}{dz}\br)^n e^z
\ee
and were first introduced in a probabilistic context in 1939 by J. Touchard \cite{T}.
The first few $T_n(z)$ are given by
\begin{eqnarray*}
T_0(z)&=&1,\quad T_1(z)=z,\quad T_2(z)=z^2+z,\\
T_3(z)&=&z^3+3z^2+z,\quad T_4(z)=z^4+6z^3+7z^2+z,\\
T_5(z)&=&z^5+10z^4+25z^3+15z^2+z,\ \ldots\ .
\end{eqnarray*}
They are also a special case of the Bell polynomials $B_n(x_1, x_2, \ldots\,,x_n)$ when all the $x_j$ ($1\leq j\leq n$) are equal, namely
\[T_n(z)=B_n(z,z, \ldots ,z)\]
and possess the generating function
\bee\label{e12}
\exp\,[z(e^t-1)]=\sum_{n=0}^\infty T_n(z)\,\frac{t^n}{n!}.
\ee

An alternative representation for the Touchard polynomials is given by
\bee\label{e13}
T_n(z)=\sum_{k=0}^n S(n,k) z^k=z^n\sum_{k=0}^n S(n,n-k) z^{-k},
\ee
where $S(n,k)$ is the Stirling number of the second kind \cite[p.~624]{DLMF}.
The second form of this representation immediately produces the expansion of $T_n(z)$ for $|z|\ra\infty$ with $n$ finite. 
With the values $S(n,n)=1$ and
\[S(n,n-1)=\bl(\!\!\begin{array}{c}n\\2\end{array}\!\!\br),\quad S(n,n-2)=\frac{3n-5}{4}\,\bl(\!\!\begin{array}{c}n\\3\end{array}\!\!\br),\quad
S(n,n-3)=\bl(\!\!\begin{array}{c}n\\4\end{array}\!\!\br)\bl(\!\!\begin{array}{c}n-2\\2\end{array}\!\!\br),
\]
we obtain
\bee\label{e14}
z^{-n}T_n(z)=1+n(n-1)\bl\{\frac{1}{2z}+\frac{(3n-5)(n-2)}{24z^2}+\frac{(n-2)^2(n-3)^2}{48z^3}+O(z^{-4})\br\}
\ee
as $|z|\ra\infty$ in the sector $-\pi\leq \arg\,z\leq\pi$.

In this note we consider the asymptotic expansion of $T_n(z)$ for large $n$ and complex values of the variable $z$
by an application of the method of steepest descents applied to a contour integral representation. In our treatment $|z|$ may be finite or allowed to be large like $O(n)$. It is sufficient to consider only $0\leq\arg\,z\leq\pi$ since, from (\ref{e13}), it is seen that 
$$T_n({\overline z})={\overline {T_n(z)}},$$ where the bar denotes the complex conjugate. We shall find that there is an infinite number of saddle points of the integrand but that the precise number that contribute to the expansion of $T_n(z)$ depends on the values of $n$ and $|z|$. In addition, when $\arg\,z=\pi$ there is a coalescence of two contributory saddle points where the neighbouring Poincar\'e-type expansions break down. Some numerical examples are given to illustrate the accuracy of the various expansions.

\vspace{0.6cm}

\begin{center}
{\bf 2. \ An integral representation and discussion of the saddle points}
\end{center}
\setcounter{section}{2}
\setcounter{equation}{0}
\renewcommand{\theequation}{\arabic{section}.\arabic{equation}}
From (\ref{e12}) we obtain the integral representation
\[
T_n(z)=\frac{n!\,e^{-z}}{2\pi i}\oint \frac{e^{ze^t}}{t^{n+1}}\,dt,
\]
where the integration path is a closed circuit described in the positive sense surrounding the origin.
Since $|\exp (ze^t)|\ra 1$ as $t\ra\infty$ in $\Re (t)<0$, it follows that when $n\geq 1$ the closed path above may be opened up into a loop, which commences at $-\infty$, encircles the origin and returns to $-\infty$.
We now consider $n\ra\infty$ with the variable $|z|=x$ either finite or large like $O(n)$. We set 
\[\mu:=\frac{n}{x},\qquad \theta:=\arg\,z,\]
where $\theta\in [0,\pi]$. Then we have
\bee\label{e21}
T_{n-1}(z)=\frac{\g(n) e^{-z}}{2\pi i} \int_{-\infty}^{(0+)} e^{n\psi(t)}dt,\qquad \psi(t)\equiv\psi(t;\mu,\theta):=\frac{e^{t+i\theta}}{\mu}-\log\,t.
\ee

Saddle points of the integrand occur when $\psi'(t)=0$; that is when
\bee\label{e22}
te^t=\mu e^{-i\theta}.
\ee
When $z$ is real and $\theta=0$, there is a saddle situated on the positive real axis given by $t_{0}=W(\mu)$, where $W$ here denotes the positive part of the principal branch of the Lambert-$W$ function; see \cite[p.~111]{DLMF}.
When $\theta=\pi$ and $0<\mu<1/e$, there are two saddles on the negative real axis given by the negative values of the Lambert function. When $\mu=1/e$, these two saddles coalesce to form a double saddle point and when $\mu>1/e$ 
the saddles move off the real axis to form a complex conjugate pair; see Section 2.1.

There is an infinite number of complex roots to (\ref{e22}) given by
\bee\label{e22a}
t+\log\,t=\log\,\mu+(2\pi k-\theta)i,
\ee
where $k$ is an integer. With $M:=(\log\,\mu)^2+(2\pi k-\theta)^2$, we find that the complex saddles $t_k$ are given approximately by
\bee\label{e23}
t_{k}\simeq \log\,\mu-\fs\log M+i(2\pi k-\theta-\arctan\phi_k),\qquad \phi_k:= \frac{2\pi k-\theta}{\log\,\mu}~.
\ee
Then for large $k$ and finite $\mu$ we see that the distribution of the complex saddles  is asymptotically described by
\[
t_{k}\simeq \log\,\mu-\log (2\pi|k|\mp\theta)+i(2\pi k-\theta\mp\fs\pi)\qquad (k\ra\pm\infty).
\] 
This last result indicates that a complex saddle occurs in horizontal strips of width $2\pi$ with the real part progressively becoming more negative as $k$ increases; see Table 1. It follows from the definition of $\psi(t)$ in (\ref{e21}) and from (\ref{e22}), (\ref{e22a}) that
\[\psi(t_k)=\frac{1}{t_k}-\log\,t_k=\frac{1}{t_k}+t_k-\log\,\mu-(2\pi k-\theta)i,\]
whence, with $\omega_k:=\arg\,t_k$,
\[\Re (\psi(t_k))=\bl(\frac{1}{|t_k|}+|t_k|\br) \cos \omega_k-\log\,\mu.\]
From this last result we can deduce that the complex saddles are increasingly subdominant as the index $|k|$ increases.

\begin{table}[th]
\caption{\footnotesize{The location of the complex saddles $t_{k}$ for different $k$ and their approximate values from (\ref{e23}) when $\mu=2$, $\theta=0$. The saddles $t_{-k}={\overline t}_{k}$, where the bar denotes the complex conjugate.}}
\begin{center}
\begin{tabular}{c|r|r}
\mcol{1}{c|}{$k$} & \mcol{1}{c|}{$t_{k}$} & \mcol{1}{c}{Approximate $t_{k}$}\\
[.05cm]\hline
&& \\[-0.2cm]
1  & $-0.83431+\ \,4.53027i$   & $-1.15078+\ \,4.82226i$\\
2  & $-1.70226+10.83981i$ & $-1.83940+11.05068i$\\
3  & $-2.15691+17.15368i$ & $-2.24402+17.31552i$\\
5  & $-2.70395+29.75450i$ & $-2.75441+29.86719i$\\
10 & $-3.42265+61.20519i$ & $-3.44738+61.27209i$\\
15 & $-3.83638+92.63559i$ & $-3.85281+92.68434i$\\
[.15cm]\hline
\end{tabular}
\end{center}
\end{table}
\vspace{0.2cm}

\noindent{2.1\ {\it Topology of the steepest descent paths}
\vspace{0.2cm}

\noindent
The paths of steepest descent, which we denote by ${\cal C}_k$, and ascent through the saddles $t_k$ are given by the paths on which
\[\Im \{\psi(t)-\psi(t_{k})\}=0.\]
The steepest descent paths terminate either in the left half-plane $\Re (t)<0$ or asymptotically approach the horizontal lines $\Im (t)=(2k+1)\pi-\theta$, $k=0, \pm 1, \pm 2,\ldots\,$ in the right half-plane. The steepest ascent paths terminate either at $t=0$ or asymptotically approach the intervening horizontal lines $\Im (t)=2k\pi-\theta$ also in the right half-plane.

A typical example of the topology of the saddles and paths of steepest descent and ascent is shown in Fig.~1 for the case $\mu=4$ and different values of $\theta$. The $t$-plane is cut along the negative real axis. When $\theta=0$, the real saddle $t_0$ is given by the Lambert function $W(4)\doteq 1.20217$, with the complex saddles $t_{\pm k}$ ($k\geq 1$) forming  conjugate pairs. The steepest descent paths through $t_{\pm 1}$ pass to infinity in $\Re (t)<0$ and so are disconnected from the remaining saddles with index $k\geq 2$. 
It then follows by Cauchy's theorem that the loop path in (\ref{e21}) commencing at $-\infty$ and encircling the origin can be deformed to pass over the three saddles $t_0$ and $t_{\pm 1}$; the remaining saddles of the infinite string are non-contributory. When $\theta=\fs\pi$, the saddle $t_0$ has moved into the lower half-plane and $t_{\pm 1}$ are no longer a conjugate pair; again the integration path
can be deformed to pass over these three saddles. 
The appearance of the fourth saddle $t_2$ occurs via a Stokes phenomenon when, at the critical value $\theta\doteq 0.76994\pi$, the steepest descent path through $t_1$ connects with the saddle $t_2$; see Fig.~1(c).
Finally, when $\theta=\pi$, there are two pairs of conjugate saddles
$t_0$, $t_1$ and $t_{-1}$, $t_2$ with the steepest descent paths through $t_2$ and $t_{-1}$ passing to infinity in $\Re (t)<0$; the integration path can now be made to pass over these four saddles. 
\begin{figure}[th]
	\begin{center}	{\tiny($a$)}\includegraphics[width=0.4\textwidth]{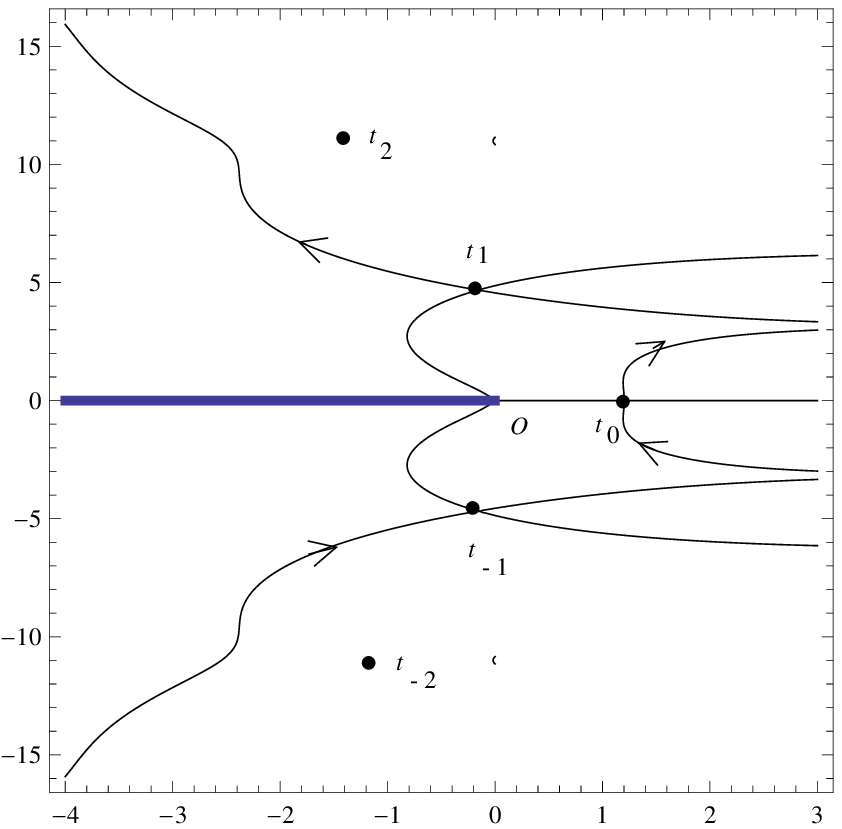}
	\qquad
	{\tiny($b$)}\includegraphics[width=0.4\textwidth]{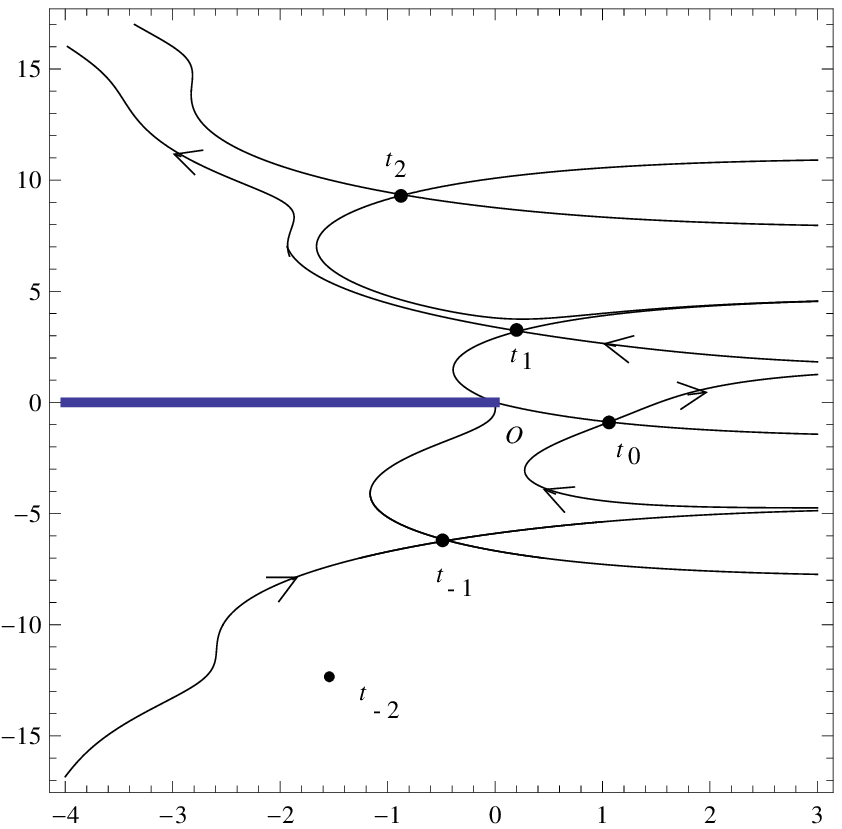}
	\vspace{0.4cm}
	
	{\tiny($c$)}\includegraphics[width=0.4\textwidth]{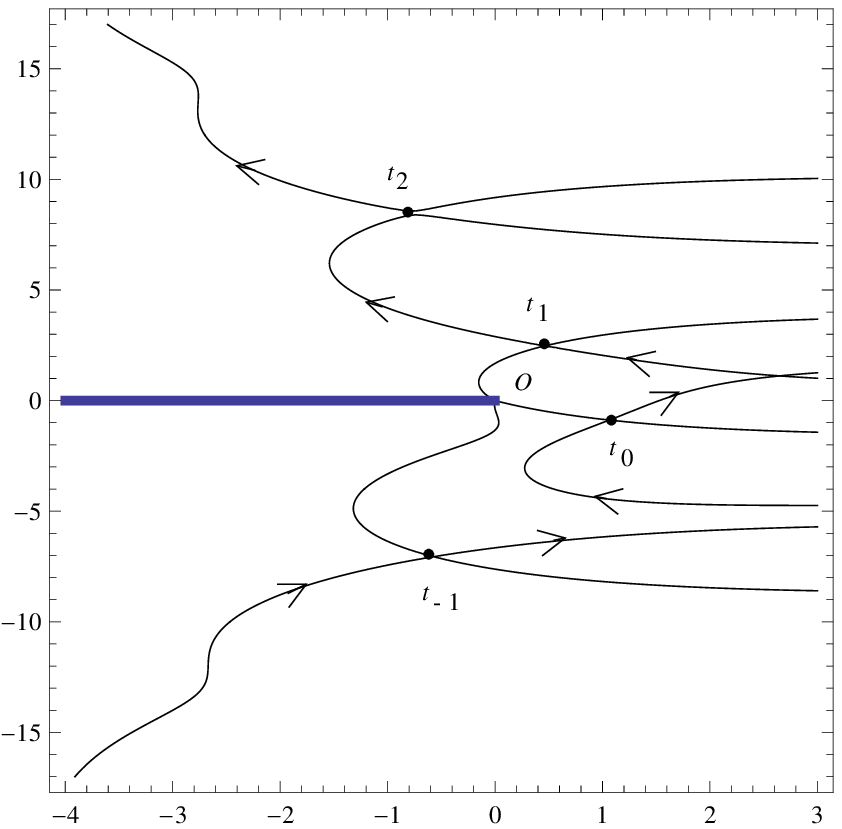}\qquad
	{\tiny($c$)}\includegraphics[width=0.4\textwidth]{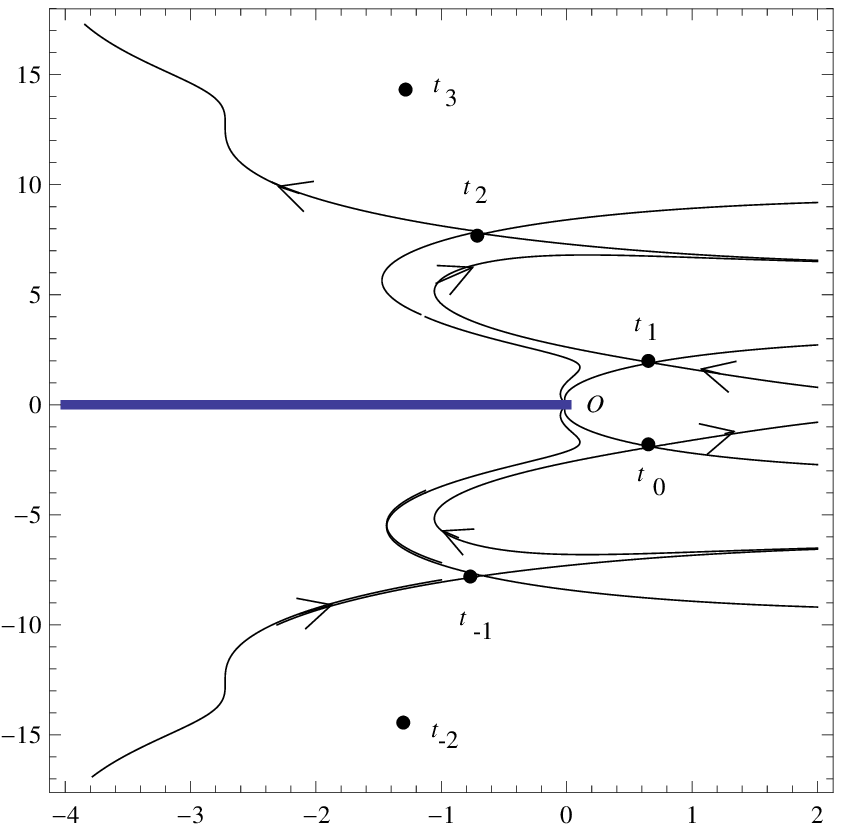}
\caption{\small{Typical paths of steepest descent and ascent through the saddles for $\mu=4$ when (a) $\theta=0$, (b) $\theta=\fs\pi$ (c) $\theta=0.76994\pi$ and (d) $\theta=\pi$. The saddles are denoted by heavy dots; the arrows indicate the direction of integration taken along steepest descent paths. There is a branch cut along $(-\infty, 0]$.
}}
	\end{center}
\end{figure}

In Fig.~2 we present an example of the steepest paths through the contributory saddles for a higher value of $\mu$.
Because of the symmetry of these paths when $\theta=0$ and $\theta=\pi$ we only show the upper half-plane; a conjugate set of paths lies in the lower half-plane. It is seen that with $\mu=12$ there are
five contributory saddles when $\theta=0$ and six saddles when $\theta=\pi$. It is also plainly visible how the steepest descent path through the last saddle ``peels away'' from the string of remaining saddles and passes to infinity in $\Re (t)<0$. 
\begin{figure}[th]
	\begin{center}	{\tiny($a$)}\includegraphics[width=0.4\textwidth]{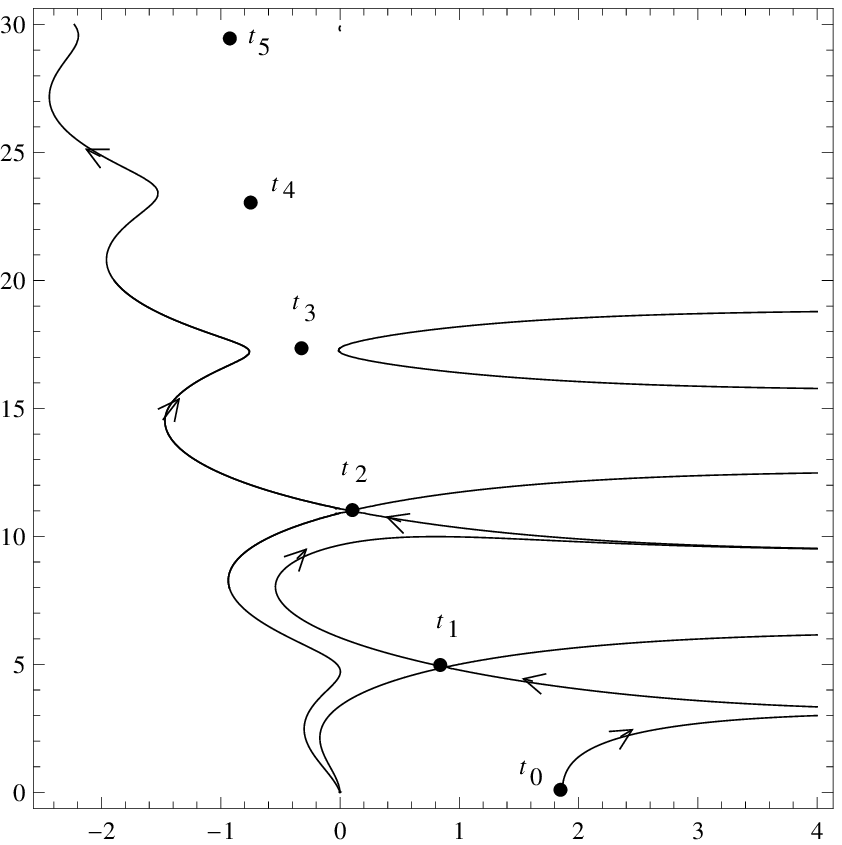}\qquad
	{\tiny($b$)}\includegraphics[width=0.4\textwidth]{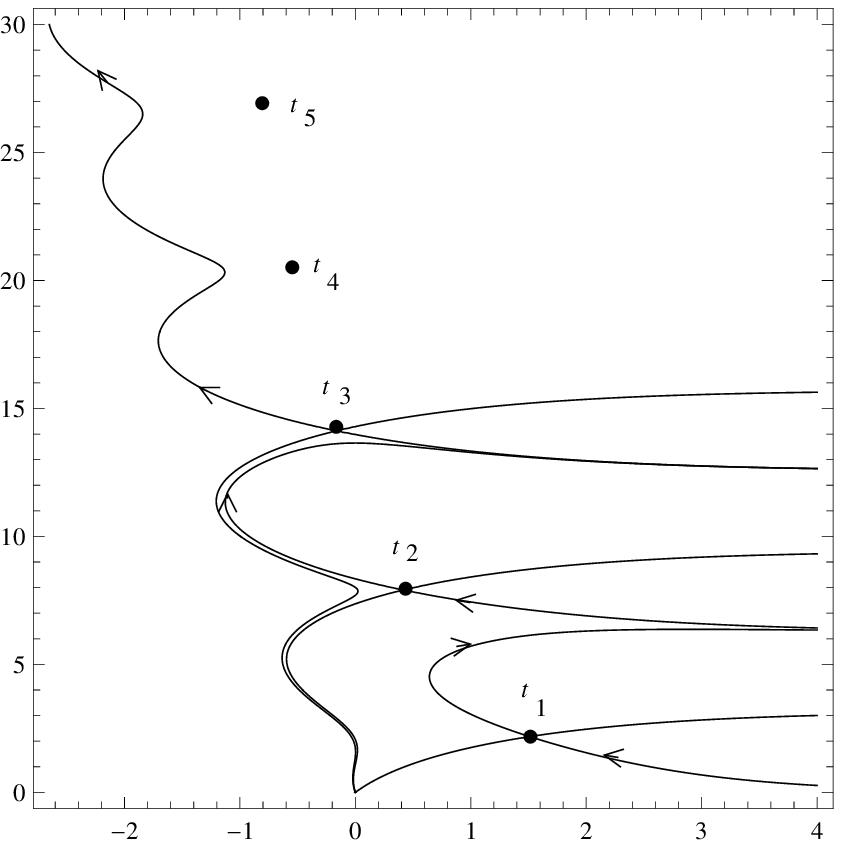}
\caption{\small{Paths of steepest descent and ascent in the upper half-plane through the contributory saddles for $\mu=12$ when (a) $\theta=0$ and (b) $\theta=\pi$; a conjugate set of paths lies in the lower half-plane. The saddles are denoted by heavy dots; the arrows indicate the direction of integration taken along steepest descent paths.}}
	\end{center}
\end{figure}
\begin{figure}[th]
	\begin{center}	{\tiny($a$)}\includegraphics[width=0.4\textwidth]{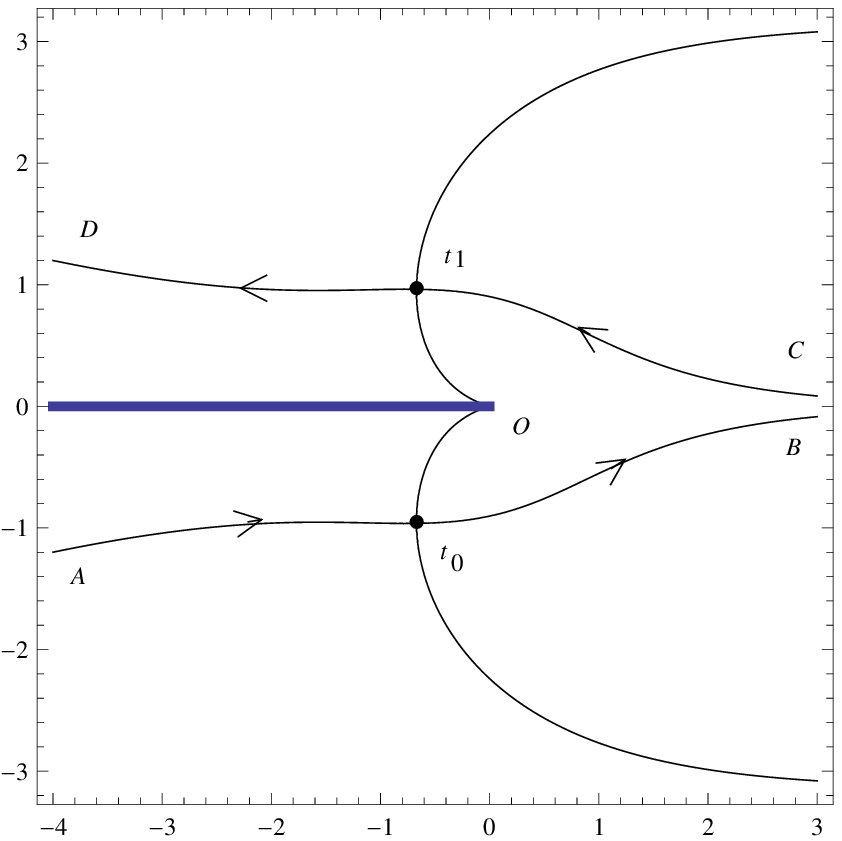}\qquad
	{\tiny($b$)}\includegraphics[width=0.4\textwidth]{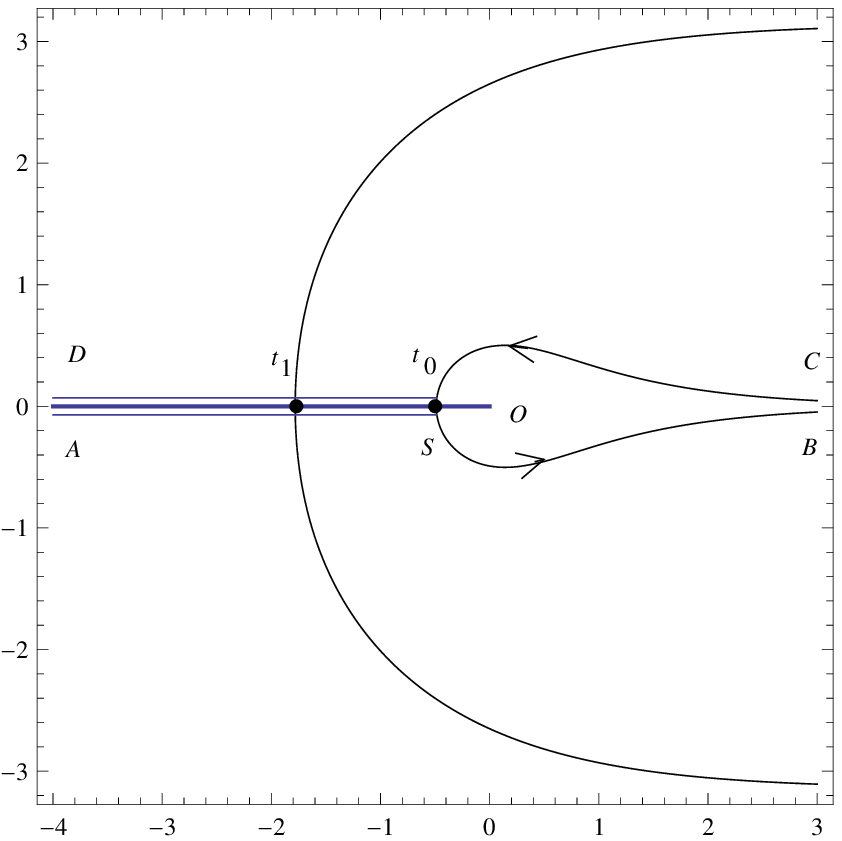}\vspace{0.4cm}
	
	{\tiny($c$)}\includegraphics[width=0.4\textwidth]{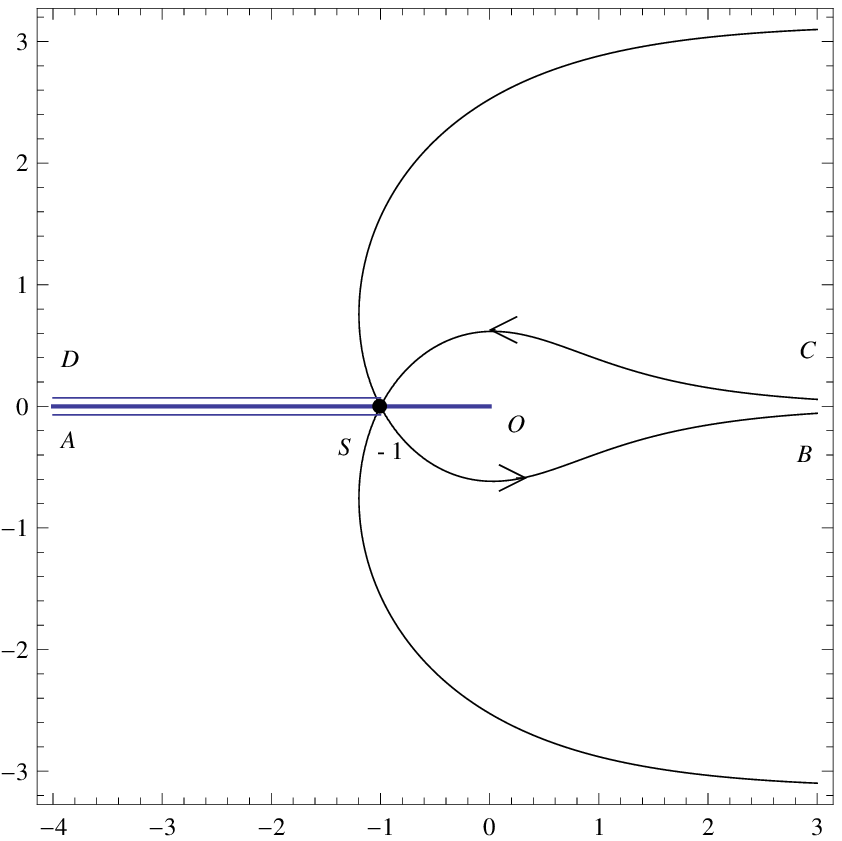}
\caption{\small{Paths of steepest descent and ascent through the saddles when (a) $1/e<\mu<\mu_1$, (b) $0<\mu<1/e$ and (c) $\mu=1/e$. The saddles are denoted by heavy dots; the arrows indicate the direction of integration  taken along steepest descent paths.  In (b) and (c) the paths $AS$ and $SD$ lie below and above the branch cut along $(-\infty, 0]$.}}
	\end{center}
\end{figure}

In Fig.~3 we show examples of the steepest paths through the contributory saddles when $\theta=\pi$ and (i) $1/e<\mu<\mu_1$, (ii) $0<\mu<1/e$ and (iii) $\mu=1/e$, where $\mu_1$ is specified below. In case (i), the saddle $t_0$ (which is on the positive axis when $\theta=0$) has rotated round the origin in the lower half-plane to form a conjugate pair with the saddle $t_1$. The integration path is the path labelled $ABCD$ in Fig.~3(a). 
When $\mu\geq\mu_1$, there are additional conjugate pairs of (subdominant) contributory saddles; see below.
In case (ii), the saddles $t_0$ and $t_1$ have rotated onto the negative real axis with $t_0\in (0,-1)$ and $t_1\in (-1,-\infty)$.
The paths of steepest descent emanating from $t_0$ pass to $+\infty$ and the paths of steepest ascent from $t_1$ asymptotically approach the lines $\Im (t)=\pm\pi$ as $\Re (t)\ra +\infty$. The integration path in (\ref{e21}) can then be deformed to pass along the lower side of the branch cut to $t_0$ and thence out to $+\infty$ along the path labelled $ASB$ in Fig.~3(b); the return path $CSD$ is the symmetrical image of that in the lower half-plane, passing to $-\infty$ along the upper side of the cut. Since $n$ is an integer, the contribution to $T_{n-1}(-x)$ from the portions of the paths along $[t_0,-\infty)$ on both sides of the cut cancel to leave\footnote{The fact that, when $z<0$, the integration path can be replaced by a loop starting and ending at $+\infty$ and encircling the origin in the positive sense can be seen from (\ref{e21}) since $|\exp (-xe^t)|\ra 0$ as $\Re (t)\ra+\infty$ when $|\Im (t)|<\fs\pi$.} the two halves of the steepest descent paths emanating from $t_0$. In case (iii), the saddles coalesce to form a double saddle at $t=-1$; the integration path then becomes the path $CSB$ in Fig.~3(c), since the contributions from $[-1,-\infty)$ on the upper and lower sides of the cut cancel.

As $\mu$ increases an increasing number of saddles contributes to the integral in (\ref{e21}). It is evident that, since $T_n(z)$ is real-valued when $\theta=0$ and $\theta=\pi$, the complex saddles must occur in conjugate pairs
with the result that there is always an odd number of contributory saddles when $\theta=0$ and an even number when $\theta=\pi$. In Tables 2 and 3 we show the number of contributory saddles when $\theta=0$ and $\theta=\pi$ for different intervals of $\mu$. The increase (or decrease) in the number of contributory saddles in a given $\mu$-interval is associated with a Stokes phenomenon that takes place at a critical value of the phase $\theta$.
Table 4 shows the values $\theta=\theta_s$ at which a Stokes phenomenon takes place as a function of $\mu$.
\begin{table}[th]
\caption{\footnotesize{Values of the interval boundaries $\mu_k$ ($1\leq k\leq 8$).}}
\begin{center}
\begin{tabular}{lr||ll}
\hline
&&&\\[-0.35cm]
\mcol{1}{c}{$k$} & \mcol{1}{c||}{$\mu_k$} & \mcol{1}{c}{$k$}  & \mcol{1}{c}{$\mu_k$}\\
[.1cm]\hline
&&&\\[-0.25cm]
1 & 3.11179  & 5 & 17.02935 \\
2 & 6.87877  & 6 & 20.13877 \\
3 & 10.25555 & 7 & 23.49898\\
4 & 13.56411 & 8 & 26.43594 \\
[.1cm]\hline
\end{tabular}
\end{center}
\end{table}
\begin{table}[th]
\caption{\footnotesize{The number of contributory saddles when $\theta=0$ and $\theta=\pi$ for different ranges of the parameter $\mu$. The end column indicates the saddle that undergoes a Stokes phenomenon.}}
\begin{center}
\begin{tabular}{c|ccc||c|crc}
\hline
&&&&&&&\\[-0.25cm]
\mcol{1}{c|}{$\mu$ interval} & \mcol{1}{c}{$\theta=0$} & \mcol{1}{c}{$\theta=\pi$} & \mcol{1}{c||}{Saddle} & \mcol{1}{c|}{$\mu$ interval}
& \mcol{1}{c}{$\theta=0$} & \mcol{1}{c}{$\theta=\pi$} & \mcol{1}{c}{Saddle}\\
[.1cm]\hline
&&&&&&&\\[-0.25cm]
$(0, \mu_1)$ & 3 & 2 & $t_{-1}$ & $(\mu_4, \mu_5)$ & 7 & 6 & $t_{-3}$\\
$(\mu_1, \mu_2)$ & 3 & 4 & $t_2$ & $(\mu_5, \mu_6)$ & 7 & 8 & $t_4$\\
$(\mu_2, \mu_3)$ & 5 & 4 & $t_{-2}$ & $(\mu_6, \mu_7)$ & 9 & 8 & $t_{-4}$\\
$(\mu_3, \mu_4)$ & 5 & 6 & $t_3$ & $(\mu_7, \mu_8)$ & 9 & 10 & $t_5$\\
[.1cm]\hline
\end{tabular}
\end{center}
\end{table}
\begin{table}[th]
\caption{\footnotesize{The values of $\theta_s$ at which a Stokes phenomenon occurs for different $\mu$.}}
\begin{center}
\begin{tabular}{ll||ll}
\hline
&&&\\[-0.25cm]
\mcol{1}{c}{$\mu$} & \mcol{1}{c||}{$\theta_s/\pi$} & \mcol{1}{c}{$\mu$} & \mcol{1}{c}{$\theta_s/\pi$} \\
[.1cm]\hline
&&&\\[-0.25cm]
0.5 & 0.26352 & 3.5 & 0.89984\\
1.0 & 0.43458 & 4.0 & 0.76994\\
1.5 & 0.57971 & 4.5 & 0.63821\\
2.0 & 0.71391 & 5.0 & 0.50446\\
3.0 & 0.97162 & 6.0 & 0.23093\\
$\mu_1$ & 1.00000 & $\mu_2$ & 0.00000\\
[.1cm]\hline
\end{tabular}
\end{center}
\end{table}

\vspace{0.6cm}

\begin{center}
{\bf 3. \ The expansion of $T_{n-1}(z)$ for large $n$ and $|z|$}
\end{center}
\setcounter{section}{3}
\setcounter{equation}{0}
\renewcommand{\theequation}{\arabic{section}.\arabic{equation}}
We denote the contribution arising from the steepest descent path ${\cal C}_k$ through the saddle $t_k$ by ${\cal J}_k(z)$, where
\[{\cal J}_k(z)=\frac{\g(n) e^{-z}}{2\pi i} \int_{{\cal C}_k} e^{n\psi(t)} dt.\]
Then since 
\[\psi(t_k)=\frac{1}{t_k}-\log\,t_k,\qquad \psi''(t_k)=\frac{1+t_k}{t_k^2},\]
application of the method of steepest descents \cite[p.~48]{DLMF} produces
\bee\label{e30}
{\cal J}_k(z)\sim\frac{\g(n) e^{-z+n/t_k}}{\sqrt{2\pi (1+t_k)}\ t_k^{n-1}}\sum_{s=0}^\infty \frac{c_{2s}(t_k) \g(s+\fs)}{n^{s+\fr} \g(\fs)}
\ee
as $n\ra\infty$. 

The coefficients $c_{2s}(t_k)$ (with $s\leq 2$) are given by \cite[p.~119]{RBD}, \cite[p.~13]{PH}
\[c_0(t_k)=1,\qquad c_2(t_k)=\frac{-1}{12\psi''(t_k)}\{5\Psi_3^2-3\Psi_4\},\]
\[c_4(t_k)=\frac{1}{864(\psi''(t_k))^2}\{385\Psi_3^4-35(6\Psi_3^2-\Psi_4)\Psi_4+168\Psi_3\Psi_5-24\Psi_6\},\]
where, for convenience in presentation, we have defined
\[\Psi_r:=\frac{\psi^{(r)}(t_k)}{\psi''(t_k)}\qquad (r\geq 3).\]
Insertion of the derivatives of $\psi(t)$ evaluated at $t_k$ then yields after some straightforward algebra the coefficients expressed in the form
\bee\label{e32}
c_2(t_k)=-\frac{P_2(t_k)}{12(1+t_k)^3},\qquad c_4(t_k)=\frac{P_4(t_k)}{864(1+t_k)^6},
\ee
where
\[P_2(t)=2t^4-3t^3-20t^2-18t+2,\]
\[P_4(t)=4t^8-156t^7-695t^6-696t^5+1092t^4+2916t^3+1972t^2-72t+4.\]

Higher-order coefficients can be obtained by an inversion process similar to that outlined in the double saddle case discussed in \cite{P}. Alternatively, they can be obtained by an expansion process to yield Wojdylo's formula \cite{W} given by
\bee\label{e320}
c_{2s}(t_k)=\frac{(-)^s}{a_0^s}\sum_{j=0}^{2s}\frac{(-)^j(s+\fs)_j}{j!\,a_0^j}\,{\cal B}_{kj}.
\ee
Here ${\cal B}_{kj}\equiv{\cal B}_{kj}(a_1, a_2, \ldots , a_{k-j+1})$ are the partial ordinary Bell polynomials generated by the recursion\footnote{For example, this generates the values $B_{41}=a_4$, $B_{42}=a_2^2+2a_1a_3$, $B_{43}=3a_1^2a_2$ and $B_{44}=a_1^4$.}
\[{\cal B}_{kj}=\sum_{r=1}^{k-j+1}a_r {\cal B}_{k-r, j-1},\qquad {\cal B}_{k0}=\delta_{k0},\]
where $\delta_{mn}$ is the Kronecker symbol and the coefficients $a_r$ appear in the expansion
\[\psi(t)-\psi(t_k)=\sum_{r=0}^\infty a_r (t-t_k)^{r+2}\]
valid in a neighbourhood of the saddle $t_k$.

\vspace{0.3cm}

\noindent{3.1\ {\it The expansion of $T_{n-1}(z)$ for $z>0$}
\vspace{0.2cm}

\noindent
From the discussion of the saddle points in Section 2.1 we then have for $x>0$
\bee\label{e33}
T_{n-1}(x)={\cal J}_0(x)+2\Re \sum_{k=1}^K {\cal J}_k(x),
\ee
where the index $K$ depends on the value of the parameter $\mu$; see Table 3. The series ${\cal J}_k(x)$ ($1\leq k\leq K$) are subdominant with respect to ${\cal J}_o(x)$ in the limit $n\ra\infty$; a correct inclusion of these contributions would necessitate the evaluation of the dominant series ${\cal J}_0(x)$ at optimal truncation (that is, truncation at, or near, the smallest term in the asymptotic series). This in turn would require the computation of the coefficients $c_{2s}(t_0)$ for large values of $s$. An example with $K=1$ is considered in Section 4.

For the moment, we neglect the subdominant contributions to yield the following result:
\newtheorem{theorem}{Theorem}
\begin{theorem}$\!\!\!.$
Let $x>0$ be either finite or at most $O(n)$. Then, neglecting exponentially smaller contributions, we have the expansion\footnote{The expansion of $T_n(x)$ is obtained from (\ref{e34}) by replacing $n$ by $n+1$.}
\bee\label{e34}
T_{n-1}(x)\sim \frac{\g(n) e^{-x+n/t_0}}{\sqrt{2\pi (1+t_0)}\ t_0^{n-1}} \sum_{s=0}^\infty \frac{c_{2s}(t_0) \g(s+\fs)}{n^{s+\fr} \g(\fs)}
\ee
as $n\ra\infty$, where $t_0=W(\mu)$ is the positive root of the equation $te^t=\mu$. The coefficients $c_{2s}(t_0)$ are specified in (\ref{e32}) for $s\leq 2$. 
\end{theorem}
A result equivalent to the leading term of (\ref{e34}) has been given for the asymptotic approximation of the probability in a Neyman type A distribution by Douglas \cite[p.~294]{D}. The relation of the leading-order approximation to the Lambert-$W$ function was pointed out by V. Vinogradov; see Remark 5.1 of \cite{VV}.

For complex $z$ with $\theta\in [0,\pi)$ the result in (\ref{e33}) is modified to
\bee\label{e35}
T_{n-1}(z)=\sum_{k=-K'}^K {\cal J}_k(z),
\ee
where the indices $K$, $K'$ depend on $\mu$ according to Table 3 and satisfy $K-K'\leq 1$. If we neglect the subdominant contributions we have 
\bee\label{e35a}
T_{n-1}(z)\sim{\cal J}_0(z)+{\cal J}_1(z)\qquad (\theta\in [0,\pi)).
\ee
For most of the $\theta$-range, ${\cal J}_1(z)$ is negligible compared to the dominant series ${\cal J}_0(z)$, except near $\theta=\pi$ where both series become comparable in importance.

\vspace{0.2cm}
\newpage
\noindent{3.2\ {\it The expansion of $T_{n-1}(z)$ for $z<0$}
\vspace{0.2cm}

\noindent
When $z<0$ ($\theta=\pi$) there are three cases to consider. First, when $1/e<\mu<\mu_1$ only the saddles $t_0$ and $t_1$, which form a conjugate pair, contribute to the integral (\ref{e21}); see Fig.~2(a). When $\mu\geq\mu_1$ there are additional pairs of conjugate saddles (see Table 3) which are subdominant as $n\ra\infty$.
It follows from (\ref{e30}) and (\ref{e35a}) that the expansion of $T_{n-1}(-x)$ when $1/e<\mu<\mu_1$ is given by
\bee\label{e36}
T_{n-1}(-x)\sim \Re \frac{\surd 2 \g(n)e^{x+n/t_0}}{\sqrt{\pi(1+t_0)}\ t_0^{n-1}} \sum_{k=0}^\infty \frac{c_{2s}(t_0) \g(s+\fs)}{n^{s+\fr} \g(\fs)}\qquad (n\ra\infty).
\ee
When $\mu\geq\mu_1$, (\ref{e36}) is the dominant expansion. 

When $0<\mu<1/e$, the saddles $t_0$ and $t_1$ are real with $t_1<t_0<0$ given by the negative roots of the Lambert-$W$ function; see Fig.~2(b). As explained in Section 2.1, only the saddle $t_0$ contributes to the integral in this case, so that
\bee\label{e37}
T_{n-1}(-x)\sim \frac{\g(n) e^{x+n/t_0}}{\sqrt{2\pi (1+t_0)}\,t_0^{n-1}}\sum_{s=0}^\infty\frac{c_{2s}(t_0) \g(s+\fs)}{n^{s+\fr} \g(\fs)}\qquad (n\ra\infty),
\ee
where $t_0$ is the smaller (negative) root of $te^t=-\mu$.

Thus we have the following theorem.
\begin{theorem}$\!\!\!.$
Let $x>0$ be either finite or at most $O(n)$. Then, we have the expansions
\[T_{n-1}(-x)\sim\left\{\begin{array}{ll}\Re \dfrac{\surd 2 \g(n)e^{x+n/t_0}}{\sqrt{\pi(1+t_0)}\ t_0^{n-1}} \sum_{s=0}^\infty \dfrac{c_{2s}(t_0) \g(s+\fs)}{n^{s+\fr} \g(\fs)} & (\mu>1/e)\\
\\
\dfrac{\g(n) e^{x+n/t_0}}{\sqrt{2\pi (1+t_0)}\,t_0^{n-1}}\sum_{s=0}^\infty\dfrac{c_{2s}(t_0) \g(s+\fs)}{n^{s+\fr} \g(\fs)} & (0<\mu<1/e) \end{array}\right.\]
as $n\ra\infty$, where $t_0$ is one of the conjugate pair of roots of $te^t=-\mu$ with smallest modulus in the first expression and the smaller (negative) root in the second expression. The upper expansion represents the dominant contribution when $\mu\geq \mu_1$.
\end{theorem}

Both the above expansions break down in the neighbourhood of $\mu=1/e$ where there is a double saddle at $t=-1$; see Fig.~2(c). A uniform asymptotic approximation valid for $\mu\sim 1/e$ and an expansion when $\mu=1/e$ as $n\ra\infty$ are discussed in \cite{P}.

\vspace{0.6cm}

\begin{center}
{\bf 4. \ Numerical examples and concluding remarks}
\end{center}
\setcounter{section}{4}
\setcounter{equation}{0}
\renewcommand{\theequation}{\arabic{section}.\arabic{equation}}
We present some numerical results to illustrate the accuracy of the expansions developed in Section 3. To keep 
the values from becoming too large we scale out the factor $n!$ and define
\bee\label{e41}
{\hat T}_n(z)=\frac{1}{n!} \sum_{k=0}^n S(n,k) z^k;
\ee
the series ${\cal J}_k(z)$ in (\ref{e30}) with a similar removal of the factor $\g(n)$ are denoted by ${\hat{\cal J}}_k(z)$.

In Table 5 we present the the values\footnote{In Tables 5--8 we have adopted the convention of writing $x(y)$ for $x\times 10^y$.} of ${\hat T}_{n-1}(x)$ and the absolute relative error in the expansion resulting from  (\ref{e34}) for different $n$ and $x$. In Table 6 we show the same when $z<0$ ($\theta=\pi$). In the first set of results with $n=20$, the values of $x\leq 50$ correspond to $\mu>1/e\doteq 0.3679$ and the expansion (\ref{e36}) applies; the remaining values $x\geq 80$ correspond to $\mu<1/e$ and so (\ref{e37}) applies. Both these expansions break down in the neighbourhood of the critical value $\mu=1/e$, which explains why the cases $n=20$, $x=50$ ($\mu=0.4$) and $n=50$, $x=150$ ($\mu=0.3$) are associated with relatively large errors. 
\begin{table}[th]
\caption{\footnotesize{Values of ${\hat T}_{n-1}(x)$ and the absolute relative error in the asymptotic expansion (\ref{e34}) for different $n$ and $x>0$ with truncation index $s=2$.}}
\begin{center}
\begin{tabular}{r|ll|ll|ll}
\hline
&&&&&&\\[-0.25cm]
\mcol{1}{c|}{} & \mcol{2}{c|}{$x=2$} & \mcol{2}{c|}{$x=5$} & \mcol{2}{c}{$x=20$}\\
\mcol{1}{c|}{$n$} & \mcol{1}{c}{${\hat T}_{n-1}(x)$} & \mcol{1}{c|}{Error} & \mcol{1}{c}{${\hat T}_{n-1}(x)$} & \mcol{1}{c|}{Error} & \mcol{1}{c}{${\hat T}_{n-1}(x)$} & \mcol{1}{c}{Error} \\
[0.1cm]\hline
&&&&&&\\[-0.25cm]
20 & 1.76101$(-02)$ & 1.713$(-05)$ & 2.07765$(+02)$ & 6.152$(-06)$ & 1.46396$(+10)$ & 3.146$(-06)$\\
30 & 2.79684$(-05)$ & 6.468$(-06)$ & 1.17615$(+01)$ & 3.119$(-06)$ & 1.10997$(+12)$ & 5.524$(-07)$\\
50 & 3.52691$(-12)$ & 1.766$(-06)$ & 8.87071$(-04)$ & 1.059$(-06)$ & 2.97967$(+13)$ & 8.697$(-08)$\\
80 & 2.20563$(-24)$ & 5.071$(-07)$ & 2,99336$(-12)$ & 3.469$(-07)$ & 2.38019$(+12)$ & 8.859$(-08)$\\
100& 2.66821$(-33)$ & 2.762$(-07)$ & 7.50809$(-19)$ & 1.987$(-07)$ & 2.04887$(+10)$ & 6.442$(-08)$\\
[.1cm]\hline
\end{tabular}
\end{center}
\end{table}
\begin{table}[th]
\caption{\footnotesize{Values of ${\hat T}_{n-1}(-x)$ and the absolute relative error in the asymptotic expansions (\ref{e36}) and (\ref{e37}) for different $n$ and $x>0$ with truncation index $s=2$.}}
\begin{center}
\begin{tabular}{r|ll|ll}
\hline
&&&&\\[-0.25cm]
\mcol{1}{c|}{} & \mcol{2}{c|}{$n=20$} & \mcol{2}{c}{$n=50$}\\
\mcol{1}{c|}{$x$} & \mcol{1}{c}{${\hat T}_{n-1}(-x)$} & \mcol{1}{c|}{Error} & \mcol{1}{c}{${\hat T}_{n-1}(-x)$} & \mcol{1}{c}{Error} \\
[0.1cm]\hline
&&&&\\[-0.25cm]
20 & $+1.72015(+03)$ & 2.144$(-04)$ & $+3.98563(-04)$ & 2.275$(-05)$\\
50 & $-1.11431(+13)$ & 1.762$(-01)$ & $+4.62064(+09)$ & 1.520$(-05)$\\
80 & $-9.07949(+17)$ & 8.459$(-04)$ & $-1.59622(+20)$ & 5.736$(-06)$\\
100& $-1.15125(+20)$ & 1.205$(-04)$ & $-1.56025(+26)$ & 6.223$(-05)$\\
150& $-5.25213(+23)$ & 1.064$(-05)$ & $-1.58180(+39)$ & 1.341$(-01)$\\
[.1cm]\hline
\end{tabular}
\end{center}
\end{table}

Table 7 presents values of ${\hat T}_{n-1}(z)$ for $n=50$ and complex $z=5e^{i\theta}$. Here the value of $\mu=10$, so that from Table 3 there are 5 contributory saddles when $\theta=0$, which reduce to 4 saddles (via a Stokes phenomenon) when $\theta=\pi$. In the asymptotic approximation we use (\ref{e35a}), which retains only the dominant series ${\hat{\cal J}}_0(z)$ together with the series ${\hat{\cal J}}_1(z)$. This latter series becomes comparable with ${\hat{\cal J}}_0(z)$ as $\theta\ra\pi$. The last column indicates the relative importance of the contribution of ${\hat{\cal J}}_1(z)$ as $\theta$ varies.
\begin{table}[th]
\caption{\footnotesize{Values of ${\hat T}_{n-1}(z)$ and its asymptotic approximation ${\hat{\cal J}}_0(z)+{\hat{\cal J}}_1(z)$
for different $\theta$ when $n=50$ and $z=5e^{i\theta}$ with truncation index $s=2$. The final column indicates the relative importance of the two asymptotic series.}}
\begin{center}
\begin{tabular}{l|l|l|c}
\hline
&&&\\[-0.25cm]
\mcol{1}{c|}{$\theta/\pi$} & \mcol{1}{c|}{${\hat T}_{n-1}(z)$} & \mcol{1}{c|}{${\hat{\cal J}}_0(z)+{\hat{\cal J}}_1(z)$} & \mcol{1}{c}{$|{\hat{\cal J}}_1(z)/{\hat{\cal J}}_0(z)|$} \\
[0.1cm]\hline
&&&\\[-0.25cm]
0.25 & $+1.42492(-04)-9.15007(-05)i$ & $+1.42492(-04)-9.15008(-05)i$ & $1.157(-29)$ \\
0.50 & $-2.13808(-07)-1.12648(-06)i$ & $-2.13810(-07)-1.12648(-06)i$ & $9.472(-22)$ \\
0.80 & $+7.58489(-12)-3.34872(-11)i$ & $+7.58481(-12)-3.34873(-11)i$ & $1.399(-09)$ \\
0.90 & $-1.49068(-13)-3.53028(-13)i$ & $-1.49070(-13)-3.53028(-13)i$ & $3.395(-05)$ \\
0.95 & $+2.67406(-14)+2.13515(-14)i$ & $+2.67407(-14)+2.13515(-14)i$ & $5.773(-03)$ \\
0.98 & $-4.92976(-15)+6.94349(-15)i$ & $-4.92974(-15)+6.94352(-15)i$ & $1.271(-01)$ \\
1.00 & $-5.42627(-15)$ & $-5.42628(-15)$ & $1.000(+00)$ \\
[.1cm]\hline
\end{tabular}
\end{center}
\end{table}
\begin{table}[th]
\caption{\footnotesize{Values of the coefficients $c_{2s}(t_0)$ for $1\leq s\leq 10$ when $\mu=4$ and $t_0=W(4)$.}}
\begin{center}
\begin{tabular}{l|c||r|c}
\hline
&&&\\[-0.25cm]
\mcol{1}{c|}{$s$} & \mcol{1}{c||}{$c_{2s}(t_0)$} & \mcol{1}{c|}{$s$} & \mcol{1}{c}{$c_{2s}(t_0)$} \\
[0.1cm]\hline
&&&\\[-0.25cm]
1 & $-3.8686291792(-01)$ & 6 & $-2.6842320622(-04)$ \\
2 & $+5.8050222467(-02)$ & 7 & $-2.9436829689(-04)$ \\
3 & $+2.3540750889(-02)$ & 8 & $+1.6066779690(-04)$ \\
4 & $-1.5978602246(-02)$ & 9 & $-4.6043216840(-05)$ \\
5 & $+4.2654445898(-03)$ &10 & $+5.7487568453(-06)$ \\
[.1cm]\hline
\end{tabular}
\end{center}
\end{table}

We provide one example to demonstrate that optimal truncation of the dominant series ${\hat{\cal J}}_0(x)$ yields an error comparable to the next subdominant series.
In the case $n=16$, $x=4$ ($\theta=0$) we have $\mu=4$, so that
\[{\hat T}_{n-1}(x)={\hat{\cal J}}_0(x)+2\Re\,{\hat{\cal J}}_1(x),\]
since from Table 3 there are no other contributory saddles for this value of $\mu$. The coefficients $c_{2s}(t_0)$, where $t_0=W(4)\doteq 1.2021679$, appearing in the expansion (\ref{e30}) have been determined by means of (\ref{e320}) for $s\leq 30$; see Table 8. The optimal truncation index was found to be $s_o=26$ for the above values of $n$ and $x$. The value of ${\hat T}_{n-1}(x)$ was computed to high precision from (\ref{e41}) and the optimally truncated dominant series, ${\hat{\cal J}}_0^{\,\mbox{{\footnotesize opt}}}(x)$, was subtracted from it to yield the value
\bee\label{e41a}
{\hat T}_{n-1}(x)-{\hat{\cal J}}_0^{\,\mbox{{\footnotesize opt}}}(x)=-1.344850\times 10^{-13}.
\ee
The value of $2\Re\,{\hat{\cal J}}_1(x)$, with $t_1\doteq -0.1573079+4.6787801i$ and with truncation index $s=2$, yields $-1.344958\times 10^{-13}$, which is close to the value in (\ref{e41a}), thereby confirming that the subdominant contribution is comparable to the error resulting from the optimally truncated dominant series.



\end{document}